\documentclass[11pt, oneside]{amsart}   
\usepackage{geometry}                		
\usepackage{graphicx}					
\usepackage{amsmath,amsthm,amsfonts,amssymb}
\usepackage{mathrsfs}
\usepackage{enumitem}
\usepackage{etoolbox}
\usepackage{xcolor}
\usepackage{url}
\usepackage{verbatim}

\apptocmd{\sloppy}{\hbadness 10000\relax}{}{}
\apptocmd{\sloppy}{\vbadness 10000\relax}{}{}
\usepackage[pdfpagelabels]{hyperref}

\newcommand{\RR}{\mathbb{R}}
\newcommand{\e}{\epsilon}

\numberwithin{equation}{section}
\theoremstyle{plain}
\newtheorem{theorem}{Theorem}[section]
\newtheorem{proposition}[theorem]{Proposition}

\newtheorem{corollary}[theorem]{Corollary}

\newtheorem{lemma}[theorem]{Lemma}

\theoremstyle{definition}
\newtheorem{remark}[theorem]{Remark}
\newtheorem{definition}[theorem]{Definition}

\newtheorem{big idea}[theorem]{The Big Idea}

\newcommand{\norm}[1]{\left\lVert#1\right\rVert}

\theoremstyle{ser}
\newtheoremstyle{ser}
{8pt}
{8pt}
{\it}
{}
{\sf}
{:}
{6mm}
{}

\newtheoremstyle{serr}
{8pt}
{8pt}
{\normalfont}
{}
{\sf}
{.}
{6mm}
{}

\theoremstyle{ser}

\theoremstyle{serr}

\title{Quantitative Stability for the Heisenberg-Pauli-Weyl inequality}
\author{Sean McCurdy}
\address{Department of Mathematical Sciences, Carnegie Mellon University, Pittsburgh, PA.}

\email{seanmccu@andrew.cmu.edu}
\author{Raghavendra Venkatraman}
\address{Department of Mathematical Sciences, Carnegie Mellon University, Pittsburgh, PA.}
\email{rvenkatr@andrew.cmu.edu}

\date{\today}						

\begin{document}
\maketitle
\begin{abstract}
    We prove a quantitative stability result for the Heisenberg-Pauli-Weyl inequality. This leads to a next, and next-to-next order correction terms in the inequality. 
\end{abstract}

\tableofcontents
\section{Introduction}

In this note, we prove a quantitative stability result for the Heisenberg-Pauli-Weyl uncertainty principle which formalizes the physical idea that a particle's position and momentum cannot both be precisely determined in any quantum state. These physical ideas were first elaborated without rigor in Heisenberg's groundbreaking 1927 paper \cite{Heisenberg27}, with rigorous mathematical formulation established later by Kennard \cite{Kennard27} and Weyl (who attributed it to Pauli) \cite{Wely28}. 

In addition to playing a fundamental role in Quantum physics, the Heisenberg-Pauli-Weyl uncertainty principle also plays an important role classical physics in signal analysis. A wave signal may be represented by a function $f$ which describes the amplitude of the wave signal as a function of time; alternatively, it may be represented through the Fourier transform $\hat{f}$ which describes how $f$ is composed of different frequencies.  In this context, the uncertainty principle describes ``limitations on the extent to which $f$ can be both time-limited and band-limited." While the importance of the uncertainty principle was not widely appreciated in signal analysis until the foundational work of Gabor \cite{Gabor46} in 1946, it appears as though it was understood in some sense by Norbert Weiner as early as 1925.

For excellent surveys with numerous additional references to the vast literature on the mathematical aspects of the uncertainty principle and related inequalities, we refer the reader to \cite{folland-sitaram,HavJor,ricoud-torresani}.

The aim of this paper is to prove a quantitative stability result for Heisenberg's inequality in Euclidean spaces $\RR^n, n \geqslant 1$. Quantitative stability results for classical inequalities in analysis and geometry have seen a burst of activity in recent years. We do not attempt a survey of the extensive literature, but refer the reader to \cite{carlensurvey,lieb,beckner}, and references therein.  However, we mention that the paper \cite{christ} addresses related questions for the Hausdorff-Young inequality by additive combinatorial techniques. By contrast, our analysis makes use of tools from the calculus of variations.  

\begin{theorem}(Heisenberg-Pauli-Weyl Inequality)\label{t:HI analyst} 
Let $u\in W^{1, 2}(\mathbb{R}^n)$ such that, additionally,
$\int_{\mathbb{R}^n}|x|^2 |u(x)|^2dx < \infty.$  Then, 
\begin{align} \label{e.HI}
\left(\int_{\mathbb{R}^n}|x|^2|u(x)|^2dx \right)\left(\int_{\mathbb{R}^n}|\nabla u(x)|^2dx \right) \geqslant \frac{n^2}{4} \left(\int_{\mathbb{R}^n}|u(x)|^2 \right)^2.
\end{align}
One has equality if and only if $u(x) \equiv c e^{\frac{1}{2\lambda}|x|^2}$ for almost every $ x \in \RR^n,$ for some $c \in \RR$ and $\lambda  < 0.$ We note that if we write, $u(x) = ce^{ \frac{1}{2\lambda}|x|^2}$, then,
\begin{align}
    \label{e:alpha control}
\lambda = -\frac{n}{2} \frac{\int_{\RR^n} |u(x)|^2\,dx }{\int_{\RR^n} |\nabla u(x)|^2\,dx} = - \frac{2}{n} \frac{\int_{\RR^n} |u(x)|^2|x|^2\,dx}{\int_{\RR^n} |u(x)|^2\,dx}. 
\end{align}
\end{theorem}
There are many proofs of this inequality, for example, see \cite[Appendix A]{Benedetto90}. 
In order to precisely formulate our results, we define the \textit{Heisenberg deficit} $\delta:$
\begin{definition}\label{d:deficit} 
For any $u \in \mathrm{Dom}(\delta) := W^{1, 2}(\mathbb{R}^n) \cap \{u \in L^2(\mathbb{R}^n) : \norm{xu}_{L^2} < \infty \},$ we define the \textit{Heisenberg deficit}, or simply deficit,
\begin{align*}
\delta(u) := \left(\int_{\mathbb{R}^n}|x|^2|u(x)|^2dx \right)\left(\int_{\mathbb{R}^n}|\nabla u(x)|^2dx \right) - \frac{n^2}{4} \left(\int_{\mathbb{R}^n}|u(x)|^2 \right)^2.
\end{align*}
\end{definition}
\begin{remark}
By Plancherel's theorem, denoting the Fourier transform of $u$ by $\hat{u},$ we note that $\delta(u) = \delta(\hat{u})$.
\end{remark}
We proceed to state our main results. The first, is a quantitative stability result for the Heisenberg inequality, which asserts that functions $f$ that have small deficit are $L^2-$close to Gaussians. Precisely, we define the \textit{extremal set}
\begin{align} \label{e.extremals}
    E := \{ c e^{-\alpha |x|^2} : c \in \RR, \alpha > 0 \}. 
\end{align}
Our main theorem is: 
\begin{theorem}(Quantitative Stability)\label{t:quant stability}
There exists a universal constant $C_1 > 0$ such that the following holds: for any  $u \in \mathrm{Dom}(\delta)$ such that $\norm{u}_{L^2} = 1,$ there exists a Gaussian $v^* = v^*(u)  \in E$ such that 
\begin{align}
\delta(u) \geqslant C_1\norm{u - v^*}^2_{L^2}.
\end{align}
\end{theorem}
The proof of Theorem \ref{t:quant stability} immediately gives rise to the following sharpening of the Heisenberg-Pauli-Weyl inequality asserted in Theorem \ref{t:HI analyst}. It contains the "next'' and ``next-to-next'' order corrections in the Heisenberg inequality. As we describe below, this inequality has been long known in one-dimension. However, to the best of our knowledge, it is new in higher dimensions. 

\begin{corollary}(Sharpened Heisenberg-Pauli-Weyl Inequality)\label{c:sharp HI}
There exists $c_4(n) > 0,$ universal, such that the following holds: for all $u \in \mathrm{Dom}(\delta)$, there exists $v^* = v^*(u) \in E$ such that the following inequality holds:
\begin{align*}
    \left(\int_{\RR^n} |x|^2 |u(x)|^2 \,dx \right)\left( \int_{\RR^n} |\nabla u(x)|^2 \,dx \right) &\geqslant \frac{n^2}{4}\left( \int_{\RR^n} |u(x)|^2\,dx\right)^2 \\ & \quad + C_1 \left(\int_{\RR^n} |u(x)|^2\,dx\right) \left( \int_{\RR^n} |u(x) - v^*(x)|^2\,dx\right)\\
    & \quad + c_4(n)\left( \int_{\RR^n} |u(x) - v^*(x)|^2\,dx\right)^2. 
\end{align*}
\end{corollary}
In order to place Corollary \ref{c:sharp HI} in context, we recall the following classical one-dimensional result of de Bruijn:

\begin{theorem}(de Bruijn, \cite{deBruijn})\label{t:deBruijn}
Let $f \in L^2(\mathbb{R}; \mathbb{C})$ and $\norm{f}_{L^2}=1$. Let $\delta > 0$ be such that for all $c>0$ and $\lambda \in \mathbb{C}$ such that $|\lambda| = 1$, we have
\begin{align*}
\norm{f - \lambda \frac{2^{\frac{1}{4}}}{c^{\frac{1}{4}}} e^{-\pi c^2t^2}}_{L^2} \geqslant \delta.
\end{align*}
Then 
\begin{align*}
\left( \int_{-\infty}^\infty t^2 |f(t)|^2\,dt\right)^{1/2} \left( \int_{-\infty}^\infty \omega^2 |\hat f(\omega)|^2\,d\omega\right)^{1/2} \geqslant \frac{1}{4\pi}[3 - 2(1 - \frac{1}{2}\delta^2)^2].
\end{align*}
\end{theorem}

The proof of Theorem \ref{t:deBruijn} proceeds by expanding $\norm{tf}_{L^2}$ and $\norm{\xi \hat{f}}_{L^2}$ in terms of the Hermite polynomials $H_{n}$, which form an orthogonal basis for $L^2(\mathbb{R}^1)$ and diagonalize the Fourier transform as an operator; see \cite[Section 1.7]{folland_HarmonicAnalysis(book)}. Using the recurrence relations satisfied by Hermite polynomials, de Bruijn shows that 
\begin{align*}
\norm{xf}_{L^2}^2 + \norm{\xi \hat{f}}_{L^2}^2 \geqslant \frac{1}{2\pi} \sum_{n = 0}^{\infty}|\langle f, H_n \rangle|^2(2n +1).
\end{align*}
Of course, there are higher-dimensional versions of Chebychev-Hermite polynomials, see, for example, \cite{grad}, which share many of the important properties of the Hermite polynomials in 1 dimension.  However, while we can use these polynomials to obtain an orthogonal basis $\{e^{\frac{-|x|^2}{4}}H^{(n)}_{i}\}_{i}$ for $L^2(\mathbb{R}^n)$, these functions do not satisfy the necessary properties with respect to the Fourier Transform. Therefore, a direct analog in higher-dimensions of de Bruijn's analysis is not possible. We view Corollary \ref{c:sharp HI} as a higher-dimensional analog of Theorem \ref{t:deBruijn}.  To the best of our knowledge, this sharpened Heisenberg Inequality, is new.

Owing to a scaling invariance of the functional $\delta$ to be discussed shortly, Theorem \ref{t:quant stability} is sharp in the following sense.

\begin{proposition} \label{t:cantimprove}
There does not exist a constant $0<C$ such that the following estimate 
\begin{align}\label{e:wrong}
    \delta(u) \geqslant C (\norm{u - v^*}^2_{L^2(\mathbb{R}^n)} + \norm{\nabla (u- v^*)}^2_{L^2(\mathbb{R}^n; \mathbb{R}^n)} + \norm{x(u-v^*)}^2_{L^2(\mathbb{R}^n)})
\end{align}
holds for all $u \in \text{Dom}(\delta)$ such that $\norm{u}_{L^2(\mathbb{R}^n)} = 1$.
\end{proposition}
 
This proposition is proved in Section \ref{s:proofs}. The proof of Theorem \ref{t:quant stability} follows by variational methods to study the functional $\delta$, and has three principle ingredients: 

\begin{enumerate}
    \item The first is a concentration compactness argument that is typical of problems with noncompact groups of symmetries. In our setting, there are two relevant invariances.  First, the group $\mathbb{R}^n$ acts on $\text{Dom}(\delta)$ be translation. Secondly, the functional $\delta$, and therefore the set $E$, is invariant under a family of rescalings. We define these rescalings, below.
\end{enumerate}

\begin{definition} \label{d:dilation}
For any $\lambda > 0,$ we define $\Phi_{\lambda}: L^2(\mathbb{R}^2) \rightarrow L^2(\mathbb{R}^n).$
\begin{align*}
\Phi_{\lambda}(f)(x) = \lambda^{\frac{n}{2}}f(\lambda x).
\end{align*}
The functional $\Phi_{\lambda}$ is linear.  Furthermore, for any $u \in \mathrm{Dom}(\delta),$ we have 
\begin{equation}\label{e.scaling}\begin{aligned}
\norm{u}_{L^2} & = \norm{\Phi_{\lambda}(u)}_{L^2}\\
\norm{\nabla u}_{L^2} & = \lambda \norm{\nabla \Phi_{\lambda}(u)}_{L^2}\\
\norm{xu}_{L^2} & = \lambda^{-1}\norm{x\Phi_{\lambda}(u)}_{L^2}.
\end{aligned}\end{equation}
In particular, we remark that $\delta(u) = \delta(\Phi_{\lambda}(u))$ and that for any $v \in E,$ $\Phi_{\lambda}(v) \in E.$ 
\end{definition}
The concentration compactness argument is contained in Section \ref{s.cc}.
\begin{enumerate}
\item[(2)] The second ingredient in the proof is a detailed study of the geometry of the extremal set $E.$ This is carried out in Section \ref{s:geometry}. In particular, using an orthogonal decomposition of a function into its radial and spherical parts, we characterize functions whose $L^2-$nearest point in $E$ is the origin, as being purely spherical functions. 
\item[(3)] The final ingredient, contained in Section \ref{s:variation}, is a precise expansion of the the deficit functional $\delta.$ This consists of an explicit computation, along with interesting cancellations arising from the minimality conditions satisfied by an $L^2-$closest Gaussian of a function.  
\end{enumerate} 
These ingredients culminate in the proof of Theorem \ref{t:quant stability} and its corollary in Section \ref{s:proofs}.


\section*{Acknowledgements:} Sean McCurdy and Raghav Venkatraman acknowledge the Center for Nonlinear Analysis where part of this work was carried out.  The research of R.V. was partially funded by the National Science Foundation under Grant No. DMS-1411646, and an AMS-Simons Travel Grant. The first author would like to thank T. Toro and Robin Neumayer for invaluable conversations which led to this project.

\section{Concentration Compactness} \label{s.cc}


The goal of this section is to study the variational problem \begin{align*}
    \min\{ \|u - v\|_{L^2} : v \in E\},
\end{align*}
for any given $u \in \mathrm{Dom}(\delta).$ We will show that each such $u$ admits a nearest point in $E$ in the $L^2$ metric. Our basic tool will be concentration compactness, see \cite{Lions}. 

\begin{theorem}\label{t:concentration compactness classic}
Let $f_k \in L^1(\mathbb{R}^n)$ be a sequence of non-negative functions such that $\norm{f_k}_{L^1(\mathbb{R}^n)} =1.$  Then, one of the following holds.
\begin{enumerate}
\item (Compactness) For all $\epsilon >0$ there exists a $R_{\epsilon}>0$ such that for all $k \in \mathbb{N}$ there exists a $y_k \in \mathbb{R}^n$ such that,
\begin{align*}
\int_{B_{R_{\epsilon}}(y_k)} f_k \geqslant 1 - \epsilon \qquad \forall k.
\end{align*}
\item (Vanishing) For all $R >0,$
\begin{align*}
\lim_{k \rightarrow \infty}\sup_{y \in \mathbb{R}^n}\int_{B_R(y)}f_k dx = 0.
\end{align*}
\item (Dichotomy) There exists an $\ell \in (0, 1)$ such that for all $ \epsilon >0$ there exists  $k_{\epsilon} \in \mathbb{N}$ and $f_k^1, f_k^2:\mathbb{R}^n \rightarrow [0, \infty)$ such that,
\begin{align*}
\lim_{k \rightarrow \infty} \mathrm{dist}\,(\mathrm{spt}(f_k^1)\, , \, \mathrm{spt}(f^2_k)) & = \infty\\
\norm{f_k - (f^1_k + f^2_k)}_{L^1} & \leqslant \epsilon\\
|\norm{f^1_k}_{L^1} - \ell| & \leqslant \epsilon\\
|\norm{f^2_k}_{L^1} - (1-\ell)| & \leqslant \epsilon
\end{align*}
for all $k \geqslant k_{\epsilon}.$ Here, $\mathrm{spt}(f)$ denotes the support of a function $f : \RR^n \to \RR.$ 
\end{enumerate}
\end{theorem}
Using this theorem, the main result of this section is:
\begin{theorem}\label{t:concentration compactness}(Concentration Compactness)
Let $M > 0.$ Let $u_i \in \mathrm{Dom}(\delta)$ such that 
\begin{align*}
&\int_{\mathbb{R}^n}|x|^2|u_i|^2dx  < M,\\
&\int_{\mathbb{R}^n}|\nabla u_i|^2dx  < M,\\
&\norm{u_i}_{L^2} = 1.
\end{align*}
Then, there is a subsequence $\{u_j\}$ and a function $u \in W^{1,2}(\mathbb{R}^n)$ such that
\begin{align*}
&\nabla u_j \rightharpoonup \nabla u \quad  \text{  in  } L^{2}(\mathbb{R}^n; \mathbb{R}^n),\\
& u_j \rightarrow u \quad  \text{  in  } L^2(\mathbb{R}^n).
\end{align*}
and $\norm{u}_{L^2} = 1.$
\end{theorem}

\begin{proof}
We argue by using Theorem \ref{t:concentration compactness classic} with the choice $f_k = |u_k|^2.$ Since $\norm{|u_k|^2}_{L^1} = 1$, we have three possibilities.  We claim that vanishing and dichotomy cannot occur, and hence, compactness must hold. \\
\textbf{Step 1.} In this step we prove that vanishing cannot occur, arguing by contradiction.  Suppose that there is a subsequence, that we continue to denote by $\{u_k\},$ such that vanishing occurs.  Then, for all $R$ large, for all $\epsilon>0$ there exists a $k_{\epsilon, R} \in \mathbb{N}$ such that for all $k \geqslant k_{\epsilon, R},$
\begin{align*}
\int_{B_R(0)}|u_k|^2dx \leqslant \epsilon.
\end{align*}
For such $u_k,$ 
\begin{align*}
\int_{\mathbb{R}^n}|x|^2|u_k|^2dx & \geqslant \int_{\mathbb{R}^n \setminus B_R(0)}|x|^2|u_k|^2dx\\
& = R^2 \int_{\mathbb{R}^n \setminus B_R(0)}|u_k|^2dx\\
& \geqslant R^2 (1-\epsilon).
\end{align*}
For $R$ large enough, this contradicts the assumption $\int_{\mathbb{R}^n}|x|^2|u_i|^2dx < M$ for all $i \in \mathbb{N}.$\\

\textbf{Step 2.} Next we show that splitting cannot occur. Once again, we assume for the sake of contradiction that that there is a subsequence, again denoted by $\{u_k\}$ such that dichotomy occurs. This means that there exists $\ell \in (0,1),$ such that for all $\e > 0,$ there exists $k_\e \in \mathbb{N}, $ and functions $f_k^1, f_k^2 : \RR^n \to [0,\infty),$ such that $\lim_{k \rightarrow \infty} \text{dist}\,(\mathrm{spt}(f_k^1)\, , \, \mathrm{spt}(f^2_k)) = \infty$, and $\||u_k|^2 - (f_k^1 + f_k^2)\|_{L^1} \leqslant \e,$ with $|\|f_k^1\| - \ell | \leqslant \e, $ and $|\|f_k^2 \| - (1-\ell)| \leqslant \e.$  For any $R$, there exists $k_{R} \in \mathbb{N}$ such that for all $k \geqslant k_{R}$, either $\text{spt}(f^1_k) \cap B_R(0) = \emptyset$ or $\text{spt}(f^2_k) \cap B_R(0) = \emptyset.$  Thus, if $\text{spt}(f^1_k) \cap B_R(0) = \emptyset$ and 
\begin{align*}
&\norm{|u_k|^2 - (f^1_k + f^2_k)}_{L^1}  \leqslant \epsilon, \,\\
&|\norm{f^1_k}_{L^1} - \ell|  \leqslant \epsilon,
\end{align*}
for $\epsilon \ll \ell$, we may infer that $\int_{\mathbb{R}^n \setminus B_R(0)} |u_k|^2dx \geqslant \ell -2\epsilon \geqslant \frac{\ell}{2}.$  Then, we have 
\begin{align*}
\int_{\mathbb{R}^n}|x|^2|u_k|^2dx & \geqslant \int_{\mathbb{R}^n \setminus B_R(0)}|x|^2|u_k|^2dx\\
& = R^2 \int_{\mathbb{R}^n \setminus B_R(0)}|u_k|^2dx\\
& \geqslant R^2 \frac{\ell}{2}.
\end{align*}
For $R$ large enough, this contradicts the assumption $\int_{\mathbb{R}^n}|x|^2|u_i|^2dx < M$ for all $i \in \mathbb{N}.$ A symmetric argument holds for $\text{supp}(f^2_k) \cap B_R(0) = \emptyset.$

\textbf{Step 3.} It follows from Steps 1 and 2 that compactness must hold. We now argue that we can take $y_k =  0$ for all $k \in \mathbb{N}.$ Note that if $y_k$ remain bounded, we may replace $R_{\epsilon} \mapsto R_{\epsilon} + |y_k|$ in order to assume that $y_k = 0.$  Suppose then, that $y_k$ is an unbounded sequence.  Then for $\epsilon = \frac{1}{2},$ for $|y_k| \geqslant R_{\frac{1}{2}} + R,$ we have 
\begin{align*}
\int_{\mathbb{R}^n \setminus B_R(0)} |u_k|^2dx \geqslant \int_{B_{R_{\frac{1}{2}}}(0)}|u|^2dx \geqslant \frac{1}{2}.
\end{align*}
Once again repeating the argument of the preceding steps, we find \begin{align*}
\int_{\mathbb{R}^n}|x|^2|u_k|^2dx & \geqslant \int_{\mathbb{R}^n \setminus B_R(0)}|x|^2|u_k|^2dx\\
& = R^2 \int_{\mathbb{R}^n \setminus B_R(0)}|u_k|^2dx\\
& \geqslant \frac{R^2}{2}.
\end{align*}
For $R$ large enough, this contradicts the assumption $\int_{\mathbb{R}^n}|x|^2|u_i|^2dx < M$ for all $i \in \mathbb{N}.$

\textbf{Step 4.} We can conclude the proof of the theorem. As $\{u_k\}$ is a bounded sequence in $W^{1,2}(\RR^n),$ by the Banach-Alaoglu theorem it follows that, upon passing to a subsequence that we do not relabel, $u_k \rightharpoonup u \in W^{1,2}(\RR^n).$ By virtue of Step 3, we know that for all $\e > 0,$ there exists $R_\e > 0$ such that for all $k \in \mathbb{N},$ $\int_{B_{R_\e}(0)} |u_k|^2 \geqslant 1 - \e.$ It follows by the Rellich-Kondrachov theorem that upon passing to a further subsequence, if necessary, $u_j \to u$ in $L^2(B_{R_\e}(0)). $ As $\int_{\RR^n \backslash B_{R_\e}(0)} |u_j|^2 < \e$ and $\|u_j\|_{L^2} = 1,$ for all $j \in \mathbb{N},$ the assertion of the theorem follows by a standard diagonalization argument as $R \to \infty.$ In particular, $\|u\|_{L^2(\RR^n)} = 1.$ 
\end{proof}


\begin{lemma}\label{l:lsc}
Let $\{u_i\}_i \in \mathrm{Dom}(\delta)$ satisfy the bounds in Theorem \ref{t:concentration compactness}. Then, if $u \in L^2$ and $\{u_i\}_i,$ a subsequence, not relabelled, such that $u_i \rightharpoonup u$ in $W^{1, 2}(\mathbb{R}^n)$ and $u_i \to u$ in $L^2(\RR^n),$ as in Theorem \ref{t:concentration compactness}, then
\begin{align}
\delta(u) \leqslant \liminf_{i \rightarrow \infty}\delta(u_i).
\end{align}
In particular, $u \in \mathrm{Dom}(\delta).$
\end{lemma}

\begin{proof}
By assumption and Theorem \ref{t:concentration compactness}, we may choose a subsequence $u_j$ such that
\begin{align*}
u_j \rightarrow u & \text{  in } L^2(\RR^n)\\
\nabla u_j \rightharpoonup \nabla u & \text{  in  } L^2(\RR^n; \mathbb{R}^n).
\end{align*}
In particular, $\|u\|_{L^2(\RR^n)} = 1.$ 

Upon possibly passing to a further subsequence that we do not relabel, $u_j \to u$ almost everywhere. By Fatou's lemma, we have that 
\begin{align*}
    \int_{\RR^n} |x|^2 |u(x)|^2 \,dx \leqslant \liminf_{j \to \infty} \int_{\RR^n} |x|^2 |u_j(x)|^2\,dx. 
\end{align*}
We also note that the $W^{1, 2}$-norm is lower semi-continuous with respect to weak convergence. Thus, 
\begin{align*}
\delta(u) & = \norm{xu}_{L^2}^2\norm{\nabla u}^2_{L^2}- \frac{n^2}{4}\norm{u}_{L^2}^4\\
 & \leqslant \liminf_{j \rightarrow \infty}\norm{xu_j}_{L^2}^2 \liminf_{j \rightarrow \infty}\norm{\nabla u_j}^2_{L^2}- \frac{n^2}{4}\lim_{j \rightarrow \infty}\norm{u_j}_{L^2}^4\\
 & \leqslant \liminf_{j \rightarrow \infty}\left(\norm{xu_j}_{L^2}^2\norm{\nabla u_j}^2_{L^2}- \frac{n^2}{4}\norm{u_j}_{L^2}^4\right)\\
 & \leqslant \liminf_{j \rightarrow \infty} \delta(u_j).
\end{align*}
\end{proof}

\begin{corollary}\label{c:strong convergence}
Let $u_i \in \mathrm{Dom}(\delta)$ satisfy $\norm{u_i}_{L^2} = 1.$  If $\delta(u_i) \rightarrow 0$, then there is a subsequence  $\{u_j\}$, a sequence $\lambda_j \in (0, \infty)$, and an extremal $v \in E$, such that $\Phi_{\lambda_j}(u_j) \rightarrow v$ in the $\norm{\cdot}_{\delta}$-norm defined by $\norm{u}_{\delta} := \norm{u}_{L^2} + \norm{\nabla u}_{L^2} + \norm{x u}_{L^2}.$
\end{corollary}

\begin{proof}
\textbf{Step 1.} We claim that we can find $M > 0,$ a subsequence $\{u_j\}$ and $\lambda_j \in (0,\infty)$ such that $\{\Phi_{\lambda_j}(u_j)\}_j$ satisfies the bounds of Theorem \ref{t:concentration compactness}. Indeed, if the whole sequence $\{u_i\}$ does not satisfy this assumption with $\lambda_i \equiv 1,$ for some $M > 0$ independent of $i,$ then there exists a subsequence $\{u_j\}$ such that either $\|\nabla u_j\|_{L^2 (\RR^n)} \to \infty$ and $\|x u_j\|_{L^2(\RR^n)} \to 0$ or vice versa. By Plancherel's theorem and the properties of the Fourier transform, it suffices to consider the first of these possibilities. Furthermore, as $\delta(u_j) \to 0$ as $j \to \infty,$ we may assume that in addition, $\frac{n}{2} \leqslant \|\nabla u_j\|_{L^2} \|x u_j\|_{L^2} < \sqrt{\frac{n^2}{4} + 1}.$

Then, setting $\lambda_j := \norm{xu_j}_{L^2},$ we compute using \eqref{e.scaling} that  
\begin{equation}
    \begin{aligned}
   & \|\Phi_{\lambda_j} (u_j)\|_{L^2} = \|u_j\|_{L^2} = 1, \\
    &\norm{\nabla \Phi_{\lambda_j}(u_j)}_{L^2} = \lambda_j \|\nabla u_j\|_{L^2} < \sqrt{\frac{n^2}{4} + 1}, \\
    & \|x  \Phi_{\lambda_j} (u_j)\| = \frac{1}{\lambda_j}\|x u_j\| = 1. 
    \end{aligned}
\end{equation}

\textbf{Step 2:} Now that we have a sequence of functions satisfying $\{\Phi_{\lambda_j}(u_j)\}_j $ satisfying the bounds in Theorem \ref{t:concentration compactness} with $M := \sqrt{\frac{n^2}{4} + 1},$ we invoke Lemma \ref{l:lsc} to obtain a further subsequence, that we do not relabel, and a function $v$ such that $\Phi_{\lambda_j}(u_j) \rightharpoonup v$ in $W^{1,2}(\mathbb{R}^n),$ and 
\begin{align*}
    0 \leqslant \delta(v) \leqslant \liminf_{j \to \infty} \delta\big( \Phi_{\lambda_j}(u_j) \big) = 0. 
\end{align*}
It follows that $v \in E.$ \\

\textbf{Step 3.} Finally, we show that in fact, $\|\nabla \big( \Phi_{\lambda_j}(u_j) - v\big)\|_{L^2} \to 0$ and $\|x\big( \Phi_{\lambda_j} (u_j) - v\big) \|_{L^2} \to 0. $ As $\{\nabla \Phi_{\lambda_j}(u_j)\}_j $ and $\{x u_j\}_j $ are $L^2$ bounded sequences, we certainly have that $\|\nabla v\|_{L^2} \leqslant \liminf_{j \to \infty} \|\nabla \Phi_{\lambda_j} (u_j)\|_{L^2}$ and $\|x v\|_{L^2} \leqslant \liminf_{j \to \infty} \|x \Phi_{\lambda_j}(u_j)\|_{L^2}.$ If either of these were a strict inequality, we obtain a contradiction to Heisenberg's inequality: indeed by Lemma \ref{l:lsc} we would have 
\begin{align*}
    \delta(v) < \liminf_{j \to \infty} \delta(\Phi_{\lambda_j}(u_j)) = \liminf_{j \to \infty} \delta(u_j) = 0.
\end{align*}
This completes the proof of the corollary. 
\end{proof}


\section{Geometry of $E$ and Closest Extremal Conditions} \label{s:geometry}

The goal of this section is to understand the geometry of the set $E$ of extremals. After first showing that each point $u \in \text{Dom}(\delta)$ has a nearest point in $E,$ we notice that $E$ is a closed cone in $L^2(\RR^n)$ that is dense in the subspace of radial square integrable functions (Lemma \ref{l:E closed}). Using these observations we characterize all points $u \in \mathrm{Dom}(\delta)$ whose nearest point in $E$ is the origin. In the latter half of this section we derive crucial transversality conditions satisfied by a nearest point projection on to $E.$

\begin{proposition} \label{p.existence}
The set $E \subset L^2(\mathbb{R}^n)$ is closed (in the $L^2(\mathbb{R}^n)$-norm).  In particular, then, for all $u \in \mathrm{Dom}(\delta)$ there exists an extremal $v^* = v^*(u)\in E$ such that
\begin{align*}
\norm{u-v^*}_{L^2} & = \min\{\norm{u - v}_{L^2} : v \in E \}.
\end{align*}
\end{proposition}

\begin{proof}
Let $u \in L^2(\mathbb{R}^n)$.  If $u \equiv 0$, then $u \in E$ and $u = v^*$.  Therefore, without loss of generality, we assume that $\norm{u}_{L^2}=1.$  Suppose that there exists a sequence $\{v_j\} \subset E$ such that $v_j \rightarrow u$ in $L^2$ as $j \to \infty$.  By Theorem \ref{t:HI analyst} we can write $v_j(x) = c_je^{-\alpha_j|x|^2}$ with $c_j \neq 0$ and $\alpha_j > 0.$ By taking subsequences, we are able to reduce to three cases:
\begin{align*}
\lim_{j \rightarrow \infty} \alpha_j  = \infty, \quad \quad 
\lim_{j \rightarrow \infty} \alpha_j  = 0 , \quad   \mbox{ and } \quad  
\lim_{j \rightarrow \infty} \alpha_j  \in (0, \infty).
\end{align*}
We claim that the first two cases can not occur. Suppose, first, for the sake of contradiction that $\lim_{j \rightarrow \infty} \alpha_j = \infty$.  Since $\norm{u}_{L^2} = 1$ and $v_j$ converges strongly to $u$ in $L^2$, there exists a constant $N \in \mathbb{N}$ such that $\norm{v_j}_{L^2} > \frac{1}{2},$ for all $j \geqslant N$.  By \eqref{e:alpha control}, since $\alpha_j \rightarrow \infty$, we must have $\norm{\nabla v_j}_{L^2} \rightarrow 0$ and $\norm{xv_j}_{L^2} \rightarrow \infty.$ Then for any $R > 0 $ we find 
\begin{align*}
    \int_{B_R(0)} |\nabla u|^2 \,dx \leqslant \liminf_{j \to \infty}\int_{B_R(0)} |\nabla v_j|^2 \,dx \leqslant \liminf_{j \to \infty} \int_{\RR^n} |\nabla v_j|^2 \,dx = 0.
\end{align*}
Letting $R \to \infty,$ it follows that $\nabla u = 0$ almost everywhere. Thus $u$ is a constant almost everywhere which is square integrable, i.e., $u = 0$ almost everywhere, contrary to $\|u\|_{L^2} = 1.$ 

The proof that the second possibility can not occur follows by Plancherel and the foregoing argument applied to the sequences $ \hat{v}_j$ converging strongly to $\hat{u}. $ We omit the details. 

It follows that we must have $\lim_{j \to \infty} \alpha_j \in (0,\infty). $  Consider the coefficient $c_n$ in the expression $v_n(x) = c_ne^{-\alpha_n|x|^2}$. Since $\lim_{n \rightarrow \infty}\alpha_n = \alpha$ exists and $c_n = c(\alpha_n) \norm{v_n}_{L^2}$ where $c(\alpha_n) = \frac{1}{\norm{e^{-\alpha_n|x|^2}}_{L^2}}$ depends continuously on $\alpha_n$, the assumption that $\alpha_n$ converges and $\norm{v_n}_{L^2} \rightarrow 1$ implies that the coefficients $c_n$ converge, as well. 

Finally, we argue that $u(x) = ce^{-\alpha|x|^2}$ with $c = \lim_{n \rightarrow \infty} c_n$ and $\alpha = \lim_{n \rightarrow \infty} \alpha_n.$  Indeed, since 
$$\langle e^{-\alpha|x|^2}, e^{-\alpha_0|x^2|} \rangle \rightarrow \norm{e^{-\alpha_0|x|^2}}^2_{L^2}$$
as $\alpha \rightarrow \alpha_0$ for any $\alpha_0 > 0,$ and $(c - c_0)e^{-\alpha |x|^2} \rightarrow 0$ in $L^2(\mathbb{R}^n)$ as $c \rightarrow c_0$ for any $c_0 \in \RR, \alpha > 0,$ we find that 
\begin{align*}
\lim_{n \rightarrow 0}\norm{ce^{-\alpha|x|^2} - c_ne^{-\alpha_n|x|^2}}^2_{L^2} & = \lim_{n \rightarrow 0} \left(1 + \norm{v_n}^2_{L^2} - 2 \langle ce^{-\alpha|x|^2} , c_ne^{-\alpha_n}|x|^2 \rangle \right)\\
& = \lim_{n \rightarrow \infty} 2 - 2 \langle ce^{-\alpha|x|^2} , c_ne^{-\alpha_n}|x|^2 +ce^{-\alpha_n}|x|^2 - ce^{-\alpha_n}|x|^2 \rangle\\
& = \lim_{n \rightarrow \infty} 2 - 2 \langle ce^{-\alpha|x|^2} , (c_n - c)e^{-\alpha_n}|x|^2\rangle - 2\langle ce^{-\alpha|x|^2} , ce^{-\alpha_n}|x|^2 \rangle\\
& = 0
\end{align*}
Thus, $u(x) = ce^{-\alpha|x|^2}$.  In particular, $u \in E$, and so $E$ is closed in $L^2(\mathbb{R}^n)$.  The continuity of the $L^2$-norm gives the desired minimality result. 
\end{proof}

\begin{remark}
We note that $v^*$ may not be unique.
\end{remark}

\begin{lemma}(Geometry of Extremals in $L^2$)\label{l:E closed}
The family of extremals $E$ defined in \eqref{e.extremals}  enjoys the following properties.
\begin{enumerate}
\item $E$ forms a closed cone in $L^2(\mathbb{R}^n)$.
\item The span of $E$ is dense in the subspace of radial functions in $L^2(\mathbb{R}^n).$ 
\end{enumerate}
\end{lemma}

\begin{proof}
$(1)$ follows from Theorem \ref{t:HI analyst} and Proposition \ref{p.existence}. To see $(2)$, we note that the span of $E$ is the same as the algebra generated by $E$.  Since the algebra generated by $\{e^{-\alpha x^2}\}_{\alpha>0} \subset C_0(\mathbb{R}^+)$ separates points on the domain $\mathbb{R}_+ = [0, \infty),$ it is dense in $C_0(\mathbb{R}_+)$ in the uniform norm-- and hence dense in $L^2(\mathbb{R}_+)$ in the $L^2$-norm--  by the Stone-Weierstrass theorem.  Therefore, the span of $E$ is dense in the collection of radial functions in $L^2(\mathbb{R}^n)$. 
\end{proof}

Because the span of $E$ is dense in the space of radial functions, we make the following orthogonal decomposition.

\begin{definition}[Orthogonal Decomposition] \label{d.orthog}
Let $u \in L^2(\mathbb{R}^n)$.  We decompose $u$ into its radial and spherical parts according to 
\begin{align*}
u = u_r + u_s,
\end{align*}
where $u_r$ is radial, defined almost everywhere by $$u_r(x) = \frac{1}{|\partial B_{|x|}(0)|}\int_{\partial B_{|x|}(0)}u d\mathcal{H}^{n-1},$$ and $u_s = u - u_r$ satisfies $\langle u_s , g \rangle_{L^2} = 0$ for all radial functions $g \in L^2(\mathbb{R}^n).$  
\end{definition}

\begin{lemma}\label{l:ortho decomp}
Let $u \in \mathrm{Dom}(\delta)$ be decomposed as $u = u_r + u_s.$  Then,
\begin{enumerate}
\item If $v^* \in \mathrm{Dom}(\delta)$ such that $\norm{u-v^*}_{L^2} = \min\{\norm{u - v}_{L^2}: v \in E\}$,
\begin{align*}
\norm{u_r-v^*}_{L^2} = \min\{\norm{u_r - v}_{L^2}: v \in E\}.
\end{align*}
\item 
\begin{align*}
 \norm{u_s}_{L^2} \leqslant \norm{u - v^*}_{L^2} \leqslant \norm{u}_{L^2}.
\end{align*}
\end{enumerate}
\end{lemma}

\begin{proof}
To see $(1),$ we note if $u = u_r + u_s$, and $g \in L^2(\mathbb{R}^n)$ is radial, then 
\begin{align} \label{e.ortog}
\norm{u - g}^2_{L^2} = \norm{u_s}^2_{L^2} + \norm{u_r - g}^2_{L^2}.    
\end{align}
Since each $v \in E$ is radial, $v^* = v^*(u)$ minimizes $\norm{u - v}^2_{L^2}$ among extremals $v \in E$ if and only if it minimizes $\norm{u_r - v}^2_{L^2},$ as  well.

To see (2), note that the first inequality is immediate from \eqref{e.ortog}, while the second follows by testing the definition of $v^*$ with the $v \equiv 0 \in E$ competitor. 
\end{proof}

\begin{proposition}\label{c:v^* = 0}
Let $u \in \mathrm{Dom}(\delta)$ be decomposed as $u = u_r+u_s$ and $v^* = v^*(u)$ be a  closest extremal as in Proposition \ref{p.existence}.  Then, if $u_r \not \equiv 0$, then $v^* \not \equiv 0$, and $v^* \equiv 0$ if and only if $u = u_s$.
\end{proposition}

\begin{proof}
Let $u_r \not \equiv 0$ be the radial part of $u$ and without loss of generality, assume that $\norm{u_r} = 1$.  By Lemma \ref{l:E closed}, the span of $E $ is dense in the space of radial functions in $L^2(\mathbb{R}^n)$.  Therefore, there exists a finite linear combination of extremals $v \in E$ such that $\norm{\sum_{i}^Na_iv_i - u_r}^2_{L^2}< \frac{1}{4}.$ By the triangle inequality, recalling that $\|u_r\|_{L^2} = 1,$ we must have that $\norm{\sum_{i=1}^N a_i v_i}_{L^2} \geqslant \frac{1}{2}.$  It follows that 
\begin{align*}
\frac{1}{4} > \norm{\sum_{i}^Na_iv_i - u_r}^2_{L^2} & = \norm{\sum_{i}^Na_iv_i}^2_{L^2} + \norm{u_r}^2_{L^2} - 2\langle \sum_{i}^Na_iv_i, u_r \rangle\\
& \geqslant \frac{1}{4} + 1 - 2\langle \sum_{i}^Na_iv_i, u_r \rangle.
\end{align*}
Hence there must be exist some $i \in \{1, \cdots, N\}$ such that $\langle v_i, u_r\rangle \not = 0.$ Fix such an $i,$ and suppose without loss of generality that $\langle v_i, u_r \rangle > 0.$ Then, for any $c \in \left( 0, 2\frac{(v_i,u_r)_{L^2}}{\norm{v_i}_{L^2}^2}\right),$ we find that $\|c v_i - u_r\|_{L^2} < \|u_r\|_{L^2}.$ But this means that the zero function can not be the element of $E$ closest to $u.$ Therefore, if $u$ has closest extremal $v^* \equiv 0$, then $u = u_s$, since $u_r \equiv 0.$ 

Conversely, if $u = u_s$, then $u$ is orthogonal to all radial functions. Lemma \ref{l:ortho decomp} then implies that $\norm{u - v^*}_{L^2} = \norm{u}_{L^2}$, and hence, $v^* \equiv 0$.
\end{proof}

\begin{corollary}\label{c:v geometry}\label{l:u-v switch sign}
Let $u \in \mathrm{Dom}(\delta) \backslash \{E\}$, and $v^* = v^*(u) \in E$ as in Proposition \ref{p.existence}.  Then, $u-v^*$ must switch sign.  In particular, $v^* - u \not\in E$. 
\end{corollary}

\begin{proof}
We note that by \eqref{e.L2}, $\langle  v^*, u-v^* \rangle = 0.$  Since $v^* \in E$ we have only two possibilities.  Either $v^* \not \equiv 0$, in which case $v^*$ does not switch sign, and hence $u-v^*$ must.  Otherwise $v^* \equiv 0$, in which case $u- v^* = u,$ so that, by Proposition \ref{c:v^* = 0}, we have $u = u_s. $ In particular, for any radial function $g \in L^2(\RR^n),$ it follows that $\int_{\RR^n} g u_s = 0.$ We may select $g = e^{-|x|^2} > 0,$ it follows that $u_s$ must change sign unless it satisfies that $u_s \equiv 0.$ The latter possibility is ruled out by the fact that $u = u_s \not\in E.$ It follows that $u - v^* = u_s - 0 = u_s $ changes sign. 
\end{proof}



We conclude this section by deriving the transversality conditions satisfied by the point $v^*(u)$ for a given point $u \in \mathrm{Dom}(\delta).$

\begin{lemma}\label{l:minimality}
Let $u \in \mathrm{Dom}(\delta)$ and let $v^*\in E$ be an extremal nearest to $u$ in the $L^2$-metric as in Proposition \ref{p.existence}. Then,
\begin{align} \label{e.L2}
\langle v^*, u - v^* \rangle_{L^2} & = 0\\ \label{e.xL2}
\langle x v^*, x(v^* - u) \rangle_{L^2} & = 0\\ \label{e.gradL2}
\langle \nabla v^*, \nabla (v^* - u) \rangle_{L^2} & = 0
\end{align}
\end{lemma}

\begin{proof}
We note that $v^*(x) = c e^{-\alpha |x|^2}$ for some $c \in \mathbb{R}$ and $\alpha > 0$.  Since $v^*$ is the closest extremal, we have that $\frac{\partial}{\partial c} || v^* - u||_{L^2}^2 = 0 =\frac{\partial}{\partial \alpha} || v^* - u||_{L^2}^2$. We calculate these partial derivatives with respect to $c$ and $\alpha:$
\begin{align*}
0 & = \lim_{h \rightarrow 0} \frac{1}{h}\left(||(1 + h)v^* - u||^2_{L^2} - ||v^* - u||^2_{L^2} \right)\\
& = \lim_{h \rightarrow 0} \frac{1}{h}\left(||v^* - u||^2_{L^2} + 2\langle h v^*, v^* - u \rangle_{L^2} + ||h v^*||^2_{L^2} - ||v^* - u||^2_{L^2} \right)\\
& = \lim_{h \rightarrow 0} \frac{1}{h}\left(2\langle h v^*, v^* - u \rangle_{L^2} + ||h v^*||^2_{L^2}  \right)\\
& = 2\langle v^*, v^* - u \rangle_{L^2},
\end{align*}
proving \eqref{e.L2}. Next, 
\begin{align*}
0 & = \lim_{h \rightarrow 0} \frac{1}{h}\left(||e^{-h|x|^2}v^* - u||^2_{L^2} - ||v^* - u||^2_{L^2} \right)\\
& = \lim_{h \rightarrow 0} \frac{1}{h}\left( \int u^2 - 2ue^{-h|x|^2}v^* + |e^{-h|x|^2}v^*|^2 - u^2 + 2uv^* -|v*|^2 dx\right)\\
& = \lim_{h \rightarrow 0} \left( \int 2u\frac{(1-e^{-h|x|^2})}{h}v^* + \frac{(e^{-2h|x|^2} - 1)}{h}|v^*|^2 dx\right)\\
& = \left( \int 2|x|^2uv^* -2|x|^2|v^*|^2 dx\right)\\
& = 2\langle x(u -v^*), xv^*\rangle_{L^2},
\end{align*}
which is \eqref{e.xL2}. Finally, we obtain \eqref{e.gradL2} by integrating by parts, the formulas $\nabla v^*(x) = -2\alpha ce^{-\alpha |x|^2}x$ and $\Delta v^* = (4 \alpha^2 |x|^2 - 2 \alpha n) v^*,$ and using \eqref{e.L2}-\eqref{e.xL2}.
\end{proof}

\begin{remark}\label{r:geometry}
From \eqref{l:minimality}, and Cauchy-Schwarz, we conclude that for any $u \in \mathrm{Dom}(\delta),$ and a corresponding nearest extremal $v^*,$ 
\begin{align*}
\norm{v^*}_{L^2} &\leqslant \norm{u}_{L^2}\\
\norm{xv^*}_{L^2} &\leqslant \norm{xu}_{L^2}\\
\norm{\nabla v^*}_{L^2} &\leqslant \norm{\nabla u}_{L^2}.
\end{align*}
\end{remark}

\begin{lemma}\label{l:large radial part} 
Let $u \in \text{Dom}(\delta)$ such that $\norm{u}_{L^2}=1$. If $u = u_r + u_s$ as in Definition \ref{d.orthog} and $\norm{u_r}_{L^2} \geqslant \frac{1}{2}$, then there exists a constant, $0< c_2(n)$, such that
\begin{align*}
    \norm{v^*(u)}_{L^2} \geqslant c_2(n).
\end{align*}
\end{lemma}
\begin{proof}
We argue by contradiction.  Let $u_i$ be a sequence of functions such that $\norm{u_i}_{L^2} = 1$ and $\norm{(u_i)_r}^2_{L^2} \geqslant \frac{1}{2}$, but suppose that $\norm{v^*(u_i)}_{L^2} \leqslant 2^{-i}.$  Arguing as in the proof of Corollary \ref{c:strong convergence}, there exists a function $u \in E$ with $\norm{u}_{L^2} = 1$ and a subsequence of rescalings, not relabelled, such that $\Phi_{\lambda_i}(u_i) \rightarrow u$ in $L^2(\RR^n)$.  In particular, arguing as before, $\norm{\Phi_{\lambda_i}(u_i) - \Phi_{\lambda_i}(v^*_i)}_{L^2} = \norm{u_i - v^*_i}$ and $\norm{\Phi_{\lambda_i}(v^*_i)}_{L^2} = \norm{v^*_i}_{L^2} \to 0,$ where $v_j^* := v^*(u_j)$ as before. It follows that $\Phi_{\lambda_j}(v_j^*) \rightarrow 0$ in $L^2$-norm. Therefore, $v^*(u) = 0$ and $\|u\|_{L^2} = 1,$ and $\|u_r\|_{L^2}^2 \geqslant \frac{1}{2}$. But $v^*(u) = 0$ implies that $u \equiv u_s,$ and this is a contradiction.
\end{proof}

\section{An expansion of $\delta$}\label{s:variation}

We begin with the calculation of $\delta(u + \epsilon \phi).$

\begin{lemma}\label{l:variation}
Let $u, \phi \in \mathrm{Dom}(\delta)$. Then, for any $\e \in \RR,$  $\delta(u + \e \phi)$ admits the following expansion: 
\begin{align}
    \delta(u + \e \phi) = \delta(u) + \e \delta^\prime(u)(\phi) + \e^2 \delta^{\prime\prime}(u) (\phi) + \e^3\delta^{'''}(u)(\phi) + \e^4 \delta(\phi), 
\end{align}
where 
\begin{equation}
    \begin{aligned}
&\delta'(u)(\phi)  := 2\left(\int_{\mathbb{R}^n}|x|^2u\phi dx\right) \left(\int_{\mathbb{R}^n}|\nabla u|^2 dx \right)  + 2\left(\int_{\mathbb{R}^n}\nabla u \cdot \nabla \phi dx \right) \left( \int_{\mathbb{R}^n}|x|^2|u|^2 dx\right)\\
& \qquad \qquad- n^2\left(\int_{\mathbb{R}^n}|u|^2dx \right)\left( \int_{\mathbb{R}^n} u \phi dx\right),
    \end{aligned}
\end{equation}
\begin{equation}
    \begin{aligned}
    \delta^{\prime\prime}(u)(\phi)&:= \left(\int_{\mathbb{R}^n}|x|^2|\phi|^2 dx\right) \left(\int_{\mathbb{R}^n}|\nabla u|^2 dx \right) + 4\left(\int_{\mathbb{R}^n}\nabla u \cdot \nabla \phi dx \right) \left( \int_{\mathbb{R}^n}|x|^2u \phi dx\right)\\
& \qquad + \left(\int_{\mathbb{R}^n}|\nabla \phi|^2 dx \right) \left( \int_{\mathbb{R}^n}|x|^2|u|^2 dx\right)- \frac{n^2}{2}\left(\int_{\mathbb{R}^n}|u|^2dx \right)\left( \int_{\mathbb{R}^n}\phi^2 dx \right)\\
& \qquad -n^2\left(\int_{\mathbb{R}^n}u \phi dx \right)^2,
    \end{aligned}
\end{equation}
and
\begin{equation}
    \begin{aligned}
    \delta^{\prime\prime\prime}(u)(\phi) & := 2 \left( \int_{\RR^n} |x|^2 u \phi dx\right) \left( \int_{\RR^n} |\nabla \phi|^2 \right) + 2 \left( \int_{\RR^n} |x|^2 \phi^2 dx\right) \left( \int_{\RR^n} \nabla \phi \cdot \nabla u \right) \\ &\qquad - n^2\left(\int_{\RR^n} u\phi dx \right)\left(\int_{\RR^n} \phi^2 dx \right).
    \end{aligned}
\end{equation}
\end{lemma}
\begin{proof} Expand $\delta(u + \epsilon \phi).$ 
\end{proof}
\begin{remark}
The notations $\delta^{\prime}(u)(\phi), \delta^{\prime\prime}(u)(\phi),$ and $\delta^{\prime\prime\prime}(u)(\phi)$ are merely suggestive of the first, second and third variations respectively; we do not invoke any differentiability properties of the functional $\delta$. 
\end{remark}

\begin{corollary}\label{l:delta expansion}
For $u\in \mathrm{Dom}(\delta) \setminus E,$ letting $v^* = v^*(u) \in E$ denote an extremal Gaussian as given by Proposition \ref{p.existence}. We have, 
\begin{align*}
\delta(u) = \norm{u - v^*}_{L^2}^2\delta''(v^*)\left(\frac{v^* - u}{\norm{v^* - u}_{L^2}}\right) + \norm{u-v^*}_{L^2}^4\delta\left(\frac{v^* - u}{\norm{v^* - u}_{L^2}}\right).
\end{align*}
\end{corollary}

\begin{proof}
We simply write $u = v^* + u - v^*,$ and invoke Lemma \ref{l:variation}, along with the orthogonality conditions from Lemma \ref{l:minimality}. 
\end{proof}


\section{Proofs of the main theorem and its corollaries} \label{s:proofs}
In this section, we finally present the proof of the main theorem. We begin with some preliminary results concerning the second variation. 
\begin{lemma}\label{l:2var positive}
Let $u \in \mathrm{Dom}(\delta) \setminus E$ be decomposed $u = u_r + u_s$ in the sense of Definition \ref{d.orthog} and let $v^*$ be an extremal Gaussian as in Proposition \ref{p.existence}. Assume that $u_r \not \equiv 0$. Then,
\begin{align*}
\delta''(v^*)(v^* - u) >0.
\end{align*}
\end{lemma}
\begin{proof}
By Lemma \ref{c:v^* = 0}, $v^* \not \equiv 0.$ Therefore, combining Lemma \ref{l:variation} and the orthogonality conditions in Lemma \ref{l:minimality}, we obtain,
\begin{align*}
\delta''(v^*)(v^* - u) & = \left(\int |x|^2 (v^* - u)^2\, dx\right)\left(\int |\nabla v^*|^2\, dx\right) - \frac{n^2}{2} \left(\int |v^*|^2\, dx\right)\left(\int(v^* - u)^2\,dx\right) \\
& \qquad + \left(\int |\nabla(v^* - u)|^2\, dx\right)\left(\int |x|^2|v^*|^2\,dx\right) .
\end{align*}
We will argue that 
\begin{align*}
\left(\int |x|^2 (v^* - u)^2\,dx\right)&\left(\int |\nabla v^*|^2\,dx\right) + \left(\int |\nabla(v^* - u)|^2\,dx\right)\left(\int |x|^2|v^*|^2\,dx\right) \\ \quad  & \geqslant  \frac{n^2}{2} \left(\int |v^*|^2\,dx\right)\left(\int(v^* - u)^2\,dx\right).
\end{align*}
Since $v^* \in E$, we have that 
$$\int |\nabla v^*|^2\,dx = \frac{n^2}{4} \frac{\left(\int|v^*|^2\,dx\right)^2}{\int |x|^2|v^*|^2dx}.
$$
Thus, we see that,
\begin{align*}
&\left(\int |x|^2 (v^* - u)^2\,dx\right)\left(\int |\nabla v^*|^2\,dx\right) + \left(\int |\nabla(v^* - u)|^2\,dx\right)\left(\int |x|^2|v^*|^2\,dx\right)\\
& \qquad = \left(\int |x|^2 (v^* - u)^2\,dx\right)\frac{n^2}{4}\left(\frac{(\int|v^*|^2dx)^2}{\int |x|^2|v^*|^2dx}\right) + \left(\int |\nabla(v^* - u)|^2\,dx\right)\left(\int |x|^2|v^*|^2\,dx\right)\\
& \qquad \geqslant n\left(\int(|v^*|^2)\,dx\right)\sqrt{\left(\int |x|^2 (v^* - u)^2\,dx\right)\left(\int |\nabla(v^* - u)|^2\,dx\right)}\\
& \qquad \geqslant  \frac{n^2}{2} \left(\int |v^*|^2\,dx\right)\left(\int(v^* - u)^2\,dx\right).
\end{align*}
where, we have used the geometric mean-arithmetic mean inequality to obtain the penultimate line and the Heisenberg inequality to obtain the last line. The last inequality is strict unless $v^* - u$ is a Gaussian. However, Corollary \ref{c:v geometry} implies that this is not the case if $u \not\in E.$ 


\end{proof}


The next lemma asserts a qualitative nondegeneracy property of the $\delta$ functional. Precisely,
\begin{lemma}\label{l:qualitative nongeneracy}
For all $\epsilon > 0$, there is a constant $c_3(n, \epsilon)>0$ such that for all $u \in \mathrm{Dom}(\delta)$ satisfying $\norm{u}_{L^2}=1$ and $\norm{u - v^*(u)}_{L^2} \geqslant \epsilon,$
\begin{align*}
\delta(u) \geqslant c_3(n, \epsilon).
\end{align*}
\end{lemma}

\begin{proof}
We argue by contradiction. Assume that for some $\e > 0,$ there exists a sequence of functions $u_i \in \mathrm{Dom}(\delta)$ with $\|u_i\|_{L^2} = 1,$ such that $\norm{u_i - v^*(u_i)}_{L^2} \geqslant \epsilon$ but for which $\delta(u_i) \leqslant 2^{-i}.$  By Corollary \ref{c:strong convergence}, there exists a Gaussian $u \in E$, a subsequence, that we do not relabel, and a sequence $\lambda_i \in (0,\infty),$ such that $\Phi_{\lambda_{i}}(u_i) \rightarrow u$ in $L^2$. However, since $\norm{\Phi_{\lambda_i}u_i - u}_{L^2} = \norm{u_i - \Phi_{\lambda_i^{-1}}(u)}_{L^2},$ and $\Phi_{\lambda_i^{-1}}(u) \in E,$ it follows that $\mathrm{dist}_{L^2(\RR^n)} (u_i, E) \to 0;$ this is contrary to our assumption that $\|u_i - v^*(u_i)\|_{L^2} \geqslant \e.$
\end{proof}

The next proposition is crucial.  It establishes a universal nondegeneracy for all functions of the form $\frac{u - v^*}{\norm{u - v^*}_{L^2}}.$

\begin{proposition}\label{t:u-v^* normalized}
Let $u \in \mathrm{Dom}(\delta) \setminus E$  with $\norm{u}_{L^2} = 1$.  Then, there exists a constant $c_4(n)>0$ such that
\begin{align*}
\delta\left(\frac{u - v^*}{\norm{u-v^*}_{L^2}}\right) > c_4(n),
\end{align*}
where $v^* = v^*(u)$ is as in Proposition \ref{p.existence}.
\end{proposition}

\begin{proof}
We argue by contradiction. Suppose that $\{u_i\} \subset \text{Dom}(\delta)$ such that $\norm{u_i}_{L^2} = 1$ and 
$$
\delta\left(\frac{u_i - v_i^*}{\norm{u_i-v_i^*}_{L^2}}\right) \leqslant 2^{-i}.
$$
Here, for brevity, we write $v^*_i := v^*(u_i)$ as granted by Proposition \ref{p.existence}. By \eqref{e.L2}, it follows that $\left\langle \frac{u_j-v_j^*}{\norm{u_j-v_j^*}_{L^2}}, v_j^* \right\rangle_{L^2} = 0$ for all $j \in \mathbb{N}.$ 

Note that by Lemma \ref{l:qualitative nongeneracy}, we may assume that $\norm{(u_i)_s}_{L^2} \le \frac{1}{2},$ since $\mathrm{dist}_{L^2}(\frac{u_i-v_i^*}{\norm{u_i - v_i^*}_{L^2}}, E) \geqslant \norm{\frac{(u_i)_s}{\norm{u_i - v_i^*}_{L^2}}} \ge \norm{(u_i)_s}_{L^2}$. Note that this implies that $\norm{(u_i)_r}_{L^2} \geqslant \frac{1}{2}$ for all sufficiently large $i \in \mathbb{N},$ and hence that $\norm{v_i^*}_{L^2} \ge c_2(n)$ by Lemma \ref{l:large radial part}.

We apply Corollary \ref{c:strong convergence} to both sequences, $\{v^*_i\}_i$ and $\left\{\frac{u_i - v_i^*}{\norm{u_i-v_i^*}_{L^2}}\right\}_i$.  Applied to the sequence $\{v^*_i\}_i$, Corollary \ref{c:strong convergence} implies that there exists a function $w_1 \in E$ such that $\norm{w_1}_{L^2}=1$ and a rescaled subsequence $\{\Phi_{\lambda_j}(v^*_j)\}_j$ such that $\Phi_{\lambda_j}(v^*_j) \rightarrow w_1$ in  $L^2(\RR^n)$. 

Applied to the sequence $\left\{\frac{u_i - v_i^*}{\norm{u_i-v_i^*}_{L^2}} \right\}_i,$ Corollary \ref{c:strong convergence} implies that there exists a function $w_2 \in E$ such that $\norm{w_2}_{L^2}=1$ and a further rescaled subsequence, still indexed $j,$ such that $\Phi_{\tau_j}\left(\frac{u_j-v_j^*}{\norm{u_j-v_j^*}_{L^2}}\right) \rightarrow w_2$ strongly in $W^{1,2}(\RR^n).$

Now, either 
\begin{align*}
\lim_{j \rightarrow \infty} \frac{\lambda_j}{\tau_j}  = \infty, \quad \quad 
\lim_{j \rightarrow \infty} \frac{\lambda_j}{\tau_j}  = 0 , \quad   \mbox{ or } \quad  
\lim_{j \rightarrow \infty} \frac{\lambda_j}{\tau_j} \in (0, \infty).
\end{align*}
In the latter case, we may apply Theorem \ref{t:concentration compactness} to the functions $\{\Phi_{\tau_j}(u_j)\}_j$ and preserve the convergence
\begin{align*}
    \Phi_{\lambda_j}\left(\frac{u_j-v_j^*}{\norm{u_j-v_j^*}_{L^2}}\right) \rightarrow \Phi_{\lambda_j}(w_2)
\end{align*}
strongly in $W^{1,2}(\mathbb{R}^n)$. By \eqref{e.L2} and the fact that $\Phi_\lambda$ is an isometry in $L^2(\mathbb{R}^n),$ we calculate,
\begin{align*}
0 & = \lim_{j \rightarrow \infty}\left\langle \frac{u_j-v_j^*}{\norm{u_j-v_j^*}_{L^2}}, v_j^* \right\rangle_{L^2}\\
& = \lim_{j \rightarrow \infty}\left\langle \Phi_{\lambda_j}\left(\frac{u_j-v_j^*}{\norm{u_j-v_j^*}_{L^2}}\right), \Phi_{\lambda_j}\left(v_j^*\right) \right\rangle_{L^2}\\
& = \langle \Phi_{\lambda_j}(w_2), w_1 \rangle_{L^2}.
\end{align*}
We obtain our contradiction by noting that $\norm{w_1}_{L^2} = \norm{\Phi_{\lambda_j}(w_2)}_{L^2} \ge c_2(n)$ implies that neither extremal is zero. 

Suppose, then, that $\lim_{j \rightarrow \infty} \frac{\lambda_j}{\tau_j}  = \infty.$  By precomposition, we may reduce to the case where $\tau_j = 1$. Thus, by \eqref{e.scaling}, we see that 
\begin{align*}
    \norm{\nabla \Phi_{\tau_j}(v_j^*)}_{L^2} \rightarrow \infty.
\end{align*}
This implies that $\Phi_{\tau_j}(v^*_j) \rightarrow \delta_0$, the Dirac mass at $0$, weakly in the sense of measures. Hence, since $\Phi_{\tau_j}\left(\frac{u_j-v_j^*}{\norm{u_j-v_j^*}_{L^2}}\right) \rightarrow w_2$ strongly in $W^{1,2}(\RR^n)$, we calculate,
\begin{align*}
0 & = \lim_{j \rightarrow \infty}\left\langle \frac{u_j-v_j^*}{\norm{u_j-v_j^*}_{L^2}}, v_j^* \right\rangle_{L^2}\\
& = \lim_{j \rightarrow \infty}\left\langle \Phi_{\tau_j}\left(\frac{u_j-v_j^*}{\norm{u_j-v_j^*}_{L^2}}\right), \Phi_{\tau_j}\left(v_j^*\right) \right\rangle_{L^2}\\
& = w_2(0) \not = 0,
\end{align*}
since $w_2$ is a nonzero Gaussian. This contradiction completes the argument in the case when $\lim_{j \to \infty} \frac{\lambda_j}{\tau_j} = \infty.$ 

The case where $\lim_{j \rightarrow \infty} \frac{\lambda_j}{\tau_j}  = 0$ is handled in an identical manner, using the symmetry under the Fourier transform. This completes the argument.  
\end{proof}

Finally, we are ready to prove the main theorem. 

\begin{proof}[Proof of Theorem \ref{t:quant stability}]
Let $u \in \mathrm{Dom}(\delta)\setminus E.$ We expand $\delta(u)$ as in Corollary \ref{l:delta expansion}, and argue as in the proof of Lemma \ref{l:2var positive} that
\begin{align*}
\delta(u) & = \norm{u-v^*}^2_{L^2}\delta''(v^*)\left(\frac{u-v^*}{\norm{u -v^*}_{L^2}}\right) + \norm{u-v^*}^4_{L^2}\delta\left(\frac{u-v^*}{\norm{u -v^*}_{L^2}}\right)\\
& \geqslant  \norm{u-v^*}^2_{L^2}\delta''(v^*)\left(\frac{u-v^*}{\norm{u -v^*}_{L^2}}\right)\\
& \geqslant \norm{u-v^*}^2_{L^2}\frac{n}{2}\norm{v^*}_{L^2}^2 \left(\sqrt{1 + \frac{4}{n^2}\delta\left(\frac{v^* - u}{\norm{v^* - u}_{L^2}}\right)} - 1\right)\\
& \geqslant \norm{u-v^*}^2_{L^2}\frac{n}{2}\norm{v^*}_{L^2}^2 \left(\sqrt{1 + \frac{4}{n^2}c_4(n)} - 1\right),
\end{align*}
where in the last line we have used Proposition \ref{t:u-v^* normalized}.   

We break the remainder of the proof into two cases using the orthogonal decomposition from Definition \ref{d.orthog}. Since $\norm{u}^2_{L^2} = \norm{u_r}^2_{L^2} + \norm{u_s}^2_{L^2} = 1$, either $\norm{u_r}^2_{L^2} \geqslant \frac{1}{2}$ or $\norm{u_r}^2_{L^2} \leqslant \frac{1}{2}$.

If $\norm{u_r}^2_{L^2} \geqslant \frac{1}{2}$, then Lemma \ref{l:large radial part} implies that there is a constant $c_2(n) > 0$ independent of $u$ such that $\norm{v^*(u)}_{L^2} \geqslant c_2(n)$.  This proves of the Theorem in the case when $\|u_r\|_{L^2}^2 \geqslant \frac{1}{2},$ with $C_1 = \frac{n}{2}c_2(n)^2(\sqrt{1 + \frac{4c_4}{n^2}} -1).$

On the other hand, if $\norm{u_r}_{L^2}^2 \leqslant \frac{1}{2},$ then $\norm{u_s}_{L^2}^2 \geqslant \frac{1}{2},$ and hence, by the orthogonality of radial and spherical functions,  $\norm{u - v^*}_{L^2}^2 \geqslant \frac{1}{2}$.  By Lemma \ref{l:qualitative nongeneracy}, there exists a constant $c_3(\frac{1}{\sqrt{2}})$ such that $\delta(u) \geqslant c_3\left(\frac{1}{\sqrt{2}}\right).$ As $\|u - v^*\|_{L^2}^2 \leqslant 1,$  it suffices to take $C_1 =c_3\left(\sqrt{\frac{1}{2}}\right)$ to obtain,
\begin{align*}
\delta(u) & \geqslant c_3\left(\frac{1}{\sqrt{2}} \right) \geqslant \norm{u-v^*}^2_{L^2} c_3\left(\frac{1}{\sqrt{2}} \right).
\end{align*}

To complete the proof of the main theorem, we let $$C_1 := \min\left(c_3\left(\sqrt{\frac{1}{2}}\right), \frac{n}{2}c_2(n)^2\left(\sqrt{1 + \frac{4c_4}{n^2}} -1\right)\right) > 0.$$ 
\end{proof}


\begin{proof}[Proof of Corollary \ref{c:sharp HI}]
Since $E$ is a cone and $v^*(cu) = c v^*(u)$, it suffices to prove Corollary \ref{c:sharp HI} for $u \in \text{Dom}(\delta)$ satisfying $\norm{u}_{L^2} =1.$ For such functions, we simply recall Corollary \ref{l:delta expansion}, Proposition \ref{t:u-v^* normalized}, Theorem \ref{t:quant stability}, and the definition of the Heisenberg deficit $\delta(u)$. The next order remainder comes directly from the expansion of $\delta$ as in Corollary \ref{l:delta expansion}. 
\end{proof}

\begin{proof}[Proof of Proposition \ref{t:cantimprove}]
We prove that for any constant given, we can produce a function which fails to satisfy \eqref{e:wrong}.  Let $C>0$ be given. Let $u \in \text{Dom}(\delta)\setminus E$ satisfying $\norm{u}_{L^2(\mathbb{R}^n)}=1.$  By the Archimedean property, for sufficiently large $\lambda > 0,$ depending upon $\norm{x(u-v^*)}_{L^2(\mathbb{R}^n )}$ and $C$,
\begin{align*}
    \delta(u) = \delta(\Phi_{\lambda}u) < \lambda C \norm{x (u-v^*)}_{L^2(\mathbb{R}^n ; \mathbb{R}^n)} = C\norm{x \Phi_{\lambda} (u-v^*)}_{L^2(\mathbb{R}^n ; \mathbb{R}^n)}.
\end{align*}
Thus, $\Phi_{\lambda}(u)$ fails to satisfy Equation $(\ref{e:wrong})$.
\end{proof}
\noindent This is surprising in light of Corollary \ref{c:strong convergence}.

\bibliographystyle{plain}
\bibliography{reference.bib}

\end{document}